Latent Markov modelling applied to grant peer review

Lutz Bornmann, Ruediger Mutz & Hans-Dieter Daniel

Dr. Lutz Bornmann

Professorship for Social Psychology and Research on Higher Education

ETH Zurich

CH-8092 Zurich

bornmann@gess.ethz.ch

Dr. Ruediger Mutz

Professorship for Social Psychology and Research on Higher Education

ETH Zurich

Prof. Dr. Hans-Dieter Daniel

Evaluation Office

University of Zurich

Professor for Social Psychology and Research on Higher Education

ETH Zurich

Abstract

In the grant peer review process we can distinguish various evaluation stages in which assessors judge applications on a rating scale. Research on this process that considers its multi-stage character scarcely exists. In this study we analyze 1954 applications for doctoral and post-doctoral fellowships of the Boehringer Ingelheim Fonds (B.I.F.), assessed in three stages (first: evaluation by an external reviewer; second: internal evaluation by a staff member; third: final decision by the B.I.F. Board of Trustees). The results show that an application only has a chance of approval if it was recommended for support in the first evaluation stage. Therefore, a form of triage or pre-screening seems desirable. We found differences in transition probabilities from one stage to the other for doctoral applicants submitted by males and females.

# 1 Introduction

An overview by the United States General Accounting Office (1999, Washington, DC, USA) of peer review practices in federal science agencies found that all of the agencies use a stepwise process with several evaluation stages, in which assessors judge applications on a rating scale. So far, however, according to our bibliographical researches, only two studies have made an empirical examination of the relationship of ratings to various stages in the grant peer review. Hodgson (1995) analyzed 779 research applications submitted to the Heart and Stroke Foundation (Ontario, Canada) from 1990 to 1994. Regression analysis established that the scores of internal reviewers (first stage) were more closely correlated to final committee score for scientific merit than those of external reviewers (second stage). Klahr (1985) analyzed nearly 200 applications that were submitted to the National Science Foundation (NSF, Arlington, Virginia, USA) and had received 1,400 reviews from "insiders" (NSF panel members) and "outsiders" (ad hoc external reviewers). The results show that ratings of the ad hoc reviewers (the external reviewers) are more "lenient" than the panel ratings. A second finding is that the outcome (approval or rejection) of about one third of the applications (the very good and the very poor) can be reliably predicted by the independent panellist assessments.

The relationship between assessors' ratings at the various evaluation stages has been determined in both studies using correlations and regression analyses. However, the time-related character of the peer review process through the successive evaluation stages is not considered by these statistical methods with several consequences: First, the estimated parameters of the simple regression analysis might be biased, because the residuals are autocorrelated due to repeated measurements. Second, information about the constancy or change of ratings for an application over multiple evaluation stages is neglected. In the various stages an application can receive the same rating (e.g., first stage: "award", second stage: "award", third stage: "award") or different ratings (e.g., first stage: "no award", second stage: "award", third stage: "possible award"). Third, measurement errors in the judgements of assessors are not taken into account.

As latent Markov models offer a good opportunity to model the peer review process adequately giving an answer to all three issues, the main aim of this study is to model quantitatively for the first time a peer review process by using latent Markov models.[1] A few years ago, the Boehringer Ingelheim Fonds (B.I.F.) agreed to have us conduct a study of its peer review process for awarding long-term doctoral and post-doctoral fellowships (Bornmann & Daniel, 2005a, 2005b, 2005c). B.I.F. is a well-known international foundation; its purpose is the promotion of basic research in biomedicine. In agreement with numerous other research institutions, the foundation uses a combination of internal and external assessments of applications in three evaluation stages for the selection of doctoral and post-doctoral fellows: first stage) evaluation by an external reviewer, second stage) internal evaluation by a staff member, and third stage) final decision by the B.I.F. Board of Trustees.

Using latent Markov models we examined a kind of test-retest reliability of the B.I.F. peer review process: the true stability of the judgements on the applications over the three evaluation stages. In addition, we analyzed initial latent class proportions and latent transition probabilities with the aim, (1) of suggesting ways of achieving a leaner peer review process, as well as (2) determining the influence of applicant's gender on the B.I.F. peer review process.

## 2   Methods

### 2.1   A multi-stage peer review process as latent Markov model

In the grant peer review process we can distinguish several evaluation stages in which both internal and external assessors judge every application on a categorical rating scale (e.g., "award", "possible award", or "no award"). Statistically speaking, every application is repeatedly measured in time on the same categorical rating scale. This kind of time series data can be well represented by Markov models, especially by latent Markov models as the most general one (Agresti, 2002; Poulsen, 1982, 1990; Singer & Spilerman, 1976/77; Wiggins, 1973). If a categorical variable representing the ratings with $J$ categories as repeated observation at, e.g., $T = 3$ evaluation stages, is available, the result is a $J^3$ contingence table. Such a table may be available not only for one group of applications but for several groups $H$ that are defined by external categorical variables (e.g., type of funding programme: doctoral and post-doctoral fellowships). According to Langeheine and van de Pol (2002) the ratings of the applications for these groups received in three evaluation stages are then given by the latent Markov model

$$p_{hijk} = \gamma_h \sum_{a=1}^{A} \sum_{b=1}^{B} \sum_{c=1}^{C} \delta_{a|h}^1 \rho_{i|ah}^1 \tau_{b|ah}^{21} \rho_{j|bh}^2 \tau_{c|bh}^{32} \rho_{k|ch}^3 \qquad (1)$$

The quantities in equation (1) are that $p_{hijk}$ is the model-expected proportion in the contingence table $(h, i, j, k)$, where $h$ denotes group $h$ (e.g., gender), and $i, j, k$ refer to the assessors' rating scale $(i = j = k = 3)$, to which the applications belong at three evaluation stages (superscripts denote the $T$ evaluation stages). Each application belongs to a group $h$, whereby the membership in $h$ remains unchanged at all $T$ evaluation stages. The proportion of group $h$ is denoted $\gamma_h$. All other parameters of the latent Markov model are conditional on group $h$.

If the same application is repeatedly assessed in different evaluation stages, it is to be expected that the assessments deviate per chance due to, e.g., attention deficits or different moods of the assessors. Human judgements are normally not fully consistent. According to Camerer and Fehr (2006) "a large body of evidence accumulated over the last three decades shows that many people violate the rationality and preference assumptions" (p. 47). With *latent* Markov models it is possible to generate ideal rating scale categories, so called latent classes, which, in contrast to the categories used by the assessors in assessing applications, (e.g., "award", "possible award", and "no award") give error-free measurement. Equation (1) shows that in the latent Markov model an application of group $h$ belongs to one of $A$ latent classes, where $\delta^I_{a|h}$ denotes the initial proportion of latent class $a$ ($a = 1, \ldots, A$). The $A$ latent classes in the Markov model are described by conditional response probabilities, where $\rho^I_{i|ah}$ is the initial probability that an application belongs to rating scale category $i$ – given membership in latent class $a$ and group $h$. If the conditional response probabilities had the value 1, the ratings of the assessors would be free of measurement error. The latent Markov model would change to a manifest Markov model.

Core elements of a latent Markov model are the transition probabilities $\tau^{2I}_{b|ah}$ (see equation 1). They allow us to quantify the proportion of applications that stay within a latent class from one evaluation step to the next and the proportion of applications that switch to another latent class. If, in a latent *stationary* Markov model, there was no switching of applications from one evaluation stage to the next – in other words the matrix of transition probabilities was an identity-matrix – the ratings would be redundant (one application would get the same assessment from the assessors in the different evaluation stages). In modelling a peer review process, however, we cannot use this stationary model, because assessors of scientific contributions frequently differ in their assessments (Cicchetti, 1991). In order to include in the modelling differences in assessment at the multiple evaluation stages, a latent *non stationary* Markov model should be estimated.

## 2.2    The multi-stage peer review process of the B.I.F.

Junior scientists submit their fellowship applications to the B.I.F. administrative office, which checks that the applicant and proposed project fulfill the formal requirements and that all required documents have been submitted (Fröhlich, 2001). Once the formal criteria have been met, the office forwards each application to an independent external reviewer – the first evaluation stage of the B.I.F. peer review process. On the basis of predetermined criteria, the reviewer assesses the application and writes a detailed review.

In the second (internal) evaluation stage of the B.I.F. peer review process a member of the foundation's staff examines the application, interviews the applicant personally and submits a detailed report. The staff member rates the application as follows: "definite award", "award", "possible award", or "no award". In the third evaluation stage the applications are submitted to the Board of Trustees. Seven internationally renowned scientists make up the Board of Trustees that convenes three times a year to make approval or rejection decisions after discussing each individual application in detail on the basis of the foregoing assessments. For approval or rejection of fellowship applications three criteria are decisive. According to Fröhlich (2001), managing director of the B.I.F., "in addition to the applicant's [track] record and the originality of the research project, there is a third element on which our judgement is based: the quality of the laboratory in which the applicant wants to pursue his project" (p. 73). For post-doctoral fellowship applicants the scientific achievements must be of outstanding quality, having resulted in papers published in or accepted by leading international journals.

## 2.3   The data set for the estimation of the latent Markov models

All in all, assessment data for 1954 applications reviewed between 1985 and 2000 were available for the calculation of the latent Markov models: 1474 applications for doctoral fellowships (75%) and 480 applications for post-doctoral fellowships (25%). The number of applications for the latter is much lower, because the foundation discontinued post-doctoral fellowships in 1995. Of the applications for doctoral fellowships 25% and of those for post-doctoral research fellowships about 20% were selected by the B.I.F. for support.

For the estimation of a latent Markov model it is assumed that for the ratings in the various evaluation stages of the peer review process categorical variables with the same rating categories are applied. As this requirement is not met in the B.I.F. peer review process (see chapter 2.2), the ratings in the various evaluation stages were commuted to a single categorical measurement system, in which two categories indicate clear decisions ("award" and "no award") and one category reflects uncertainty in reaching a decision ("possible award").

Since the external reviewers in the first evaluation stage did not use a rating scale, two experts of the Centre for Research on Higher Education and Work (Kassel, Germany) independently rated all final statements of the reviewers afterwards according to the proposed scale. The reliability of the two experts' ratings is very high (weighted kappa coefficient = .96). The four rating categories used by the staff members in the second evaluation stage of the B.I.F. peer review process are transformed into three categories, by merging "definite award" and "award" into the category "award". For the transformation of the final decisions of the Board of Trustees (approval or rejection) into a variable with three categories we proceeded as follows. At each of the three Board meetings per year, the seven members of the Board decide on applications in three rounds. In the first round of decision-making, some fellowship applications are approved (rated 'A'), some are rejected (rated 'A-B' and lower), and some are earmarked for consideration in the next round (rated 'A-'). In the second and, if necessary, third decision round, the number of applications approved or dismissed depends on how much funding is still available for the session (Fröhlich, 2001, p. 76).

Those applications earmarked in the first round for consideration in the next round just failed to persuade the Trustees (otherwise they would have been accepted immediately), but the applications were considered sufficiently promising that they were not immediately rejected. In the first round the Trustees thus use, for the rating of applications, the categories "approval", "rejection", and "decision adjourned" that could be used for the calculation of the latent Markov models (we categorized a decision to adjourn an application to the next round as "possible award").

Table 1 shows the recoded data for the B.I.F. peer review process as it is used for the calculation of the latent Markov models.

## 2.4 The likelihood-ratio test with bootstrapping

To obtain comparisons of different Markov models as they are adapted to the B.I.F. data, we used the likelihood-ratio test. Langeheine et al. (1996) offer a bootstrapping method to get a valid goodness of fit statistic in the case of a sparse data set (see Table 1), when model-expected frequencies are 0 or when model probabilities are estimated 0 or 1. The population probabilities in all cells of Table 1 (here $3 \times 3 \times 3$) can be written in one vector $P$. After sampling applications from each cell, we get a sample estimate of $P$, written $p$. A Markov model describes these proportions in terms of a limited number of parameters $\Phi$, $P = f(\Phi)$, which can be estimated by the sample $p$ with maximum likelihood (ML) algorithm. The estimated parameters $\Phi'$ allow us to derive an ML estimate of population probabilities $P' = f(\Phi')$, given sample proportion $p$ and the model. The estimated $P'$ will not be exactly equal to the true probability $P$, because a sample is drawn.

If samples are drawn from the population several times, defined by $P'$, a distribution of loglikelihood-ratios can be obtained. If the upper 5% of this distribution contains large LR values that are sparse, one can reject the hypothesis that the model really holds in the population. "Loglikelihood theory shows that in case of simple random sampling and sample size going to infinity, a $\chi^2$-distribution applies with the model's degrees of freedom as its means, $E(G^2) = df$" (Langeheine et al., 1996, pp. 494-495). A statistical program package, called PANMARK (PANel analysis using MARKov chains, van de Pol et al., 2000), allows the estimation of latent Markov models by offering bootstrapping to get valid likelihood-ratio test statistics.

# 3   Results

## 3.1   The latent non stationary Markov model

For the B.I.F. peer review process we have estimated a latent non stationary Markov model, in which the following assumptions are associated with the process (see also chapter 2):

(1)  The assessors (i.e. the independent external reviewers, the members of the foundation's staff, the Board of Trustees) have as a general rule two different types of application (doctoral and post-doctoral) to assess. Therefore, for the model estimation we have two $J^3$ contingence tables with ratings on fellowship applications.

(2)  The independent external reviewers, the members of the foundation's staff and the Board of Trustees do not assess with full consistency; they infringe in their assessments the rationality and preference assumptions (see Camerer & Fehr, 2006). For this reason we have provided, for our model estimation, latent classes with ideal rating categories. The model assumes that the amount of measurement errors in the ratings of the assessors is equal for each evaluation stage of the peer review process.

(3)  One application may be differently assessed by the assessors in the three evaluation stages (1: external reviewer, 2: staff member, 3: Board of Trustees). In the model estimation we have therefore taken the non stationary model of the ratings. This means the estimation also assumes that applications for doctoral and post-doctoral fellowships do not differ in latent class proportions and response probabilities at the beginning of the assessment process.

Other Markov models we tested with the B.I.F. peer review data and which involve other assumptions for it (e.g., a latent *stationary* Markov model or a *manifest* Markov model *without* latent classes), failed to obtain the fit of the latent non stationary Markov models (tested with the likelihood-ratio test, see chapter 2.4).

## 3.2    Reliability of the B.I.F. peer review process

Before we discuss, in the chapters 3.3 and 3.4, the parameter estimations of the latent non stationary Markov models, in this chapter we should like to present our findings on the reliability of the B.I.F. peer review process – stability and change of ratings over the three evaluation stages (see Langeheine & van de Pol, 1990). As in latent Markov models measurement errors are taken into account, we are able to decide between true change and error as well as true stability and error – similar to structural equation modelling (SEM). In the context of classical test theory in psychometry the reliability of a measure is defined as the proportion of true variability to total variability (Steyer, 1989). If we transfer this definition, the reliability of the peer review process (a kind of test-retest reliability) is given by the true proportion of those applications (1 - proportion of measurement error of change) for which the ratings received in the first evaluation stage (the external reviewer) does not change in the second (staff member) and third evaluation stage (Board of Trustees). The values of the reliability coefficient can theoretically vary between 0.00 and 1.00 and give the "true" agreement in the ratings between the B.I.F. assessors (external reviewer, staff member, Board of Trustees). If the value of the reliability coefficient of a peer review process is 1, the assessors would be to a large extent able reliably to assign a rating category to an application.

To obtain the reliability for the B.I.F. peer review process, we must take a closer look at the quantities of the full $27 \times 27$ cross table of expected frequencies (where the rows correspond to the manifest response pattern and the columns correspond to the latent class pattern). The columns 111, 222, 333 of this table (1 = "award", 2 = "possible award", 3 = "no award") obtain the total "no change" part (i.e. at each evaluation stage the assessment was the same), which can break down into stability and error components. The true stability is captured by the cells with the same row pattern (111, 222, 333); the error is defined as total stability (column sum) minus true stability. The same calculations can be made for the columns of the table giving us the "change" part (Langeheine, 1988). A macro-program with the statistic software SAS (Version 8.0) was created by the second author of this paper to calculate the reliability.

Table 2 gives an impression of the reliability of the B.I.F. peer review process separated into doctoral and post-doctoral applications. In manifest data the proportion of stability is 24% for applications for doctoral fellowships and 22% for applications for post-doctoral fellowships; the proportions of change amount correspondingly to 76% and 78%. Using the latent Markov model estimations we calculated, as described above, the measurement errors for stability and change of the process. As shown by Table 2, in both application groups the measurement error proportion for stability is very low (3% for doctoral and 2% for post-doctoral applications). The measurement errors for change in both groups are approximately 20%. In view of the fact that the reliability (e.g., Cronbach Alpha) of most of the psychological test inventories lies between 0.80 and 0.90 (Peterson, 1994), with a value of $\sim 0.80 = (1 - 0.20)$ the reliability of the B.I.F. peer review process can be regarded as sufficient.

## 3.3 Parameter estimations of the latent Markov model for the B.I.F. peer review process

The results of the latent Markov models are shown in Table 3a and Table 3b. In Table 3a we show those initial latent class proportions of doctoral and post-doctoral applications that were estimated in the Markov model as initial proportions for the transition probabilities in Table 3b (see the descriptions in the next section).

*Initial latent class proportion and response probabilities*

Table 3a shows (1) the proportions of doctoral and post-doctoral applications for three latent classes and (2) the response probabilities ("award", "possible award", and "no award") for doctoral and post-doctoral applications of belonging to one of the three latent classes at the beginning of the assessment process. As a stated assumption of the estimated latent Markov model the applications for doctoral and post-doctoral fellowships do not differ in latent class proportions and response probabilities in Table 3a. About 60% of the applications falls in the first latent class. The response probabilities (row percent) show that in this latent class the majority are applications that received an "award" rating ($\rho = 0.92$) and hardly any that received a "possible award" ($\rho = 0.06$) or "no award" rating ($\rho = 0.02$). The second latent class (18% of the applications) represents mainly applications with the initial rating "possible award" ($\rho = 0.81$); including those applications that are rated at the beginning of the assessment process as "award" ($\rho = 0.17$). The third latent class (20% of the applications) represents applications that are rated with "no award"; the response probability has the maximal value ($\rho = 1.00$).

These results suggest that the first latent class represents mainly applications rated at the beginning of the assessment process with "award", the second latent class represents applications that are rated with "possible award" and the third latent class represents those rated with "no award". The reliability of identifying these latent classes using response probabilities is very high for the first ($\rho = 0.92$) and third ($\rho = 1.00$) latent classes and moderate for the second latent class ($\rho = 0.81$).

*Latent transition probabilities*

The transition probabilities in Table 3b reveal the probabilities with the B.I.F. peer review process, when moving from one evaluation stage to the next, of remaining in the same latent class (rating category) or changing to a different one. On the left side of the table are the latent transition probabilities for the transition from the first to the second evaluation stage ($t_1$ - $t_2$, separated for doctoral and post-doctoral applications) and on the right side those for the transition from the second to the third stage ($t_2 - t_3$). So the four diagonals in the table rule off the probabilities of remaining in one latent class when making the transition from one evaluation stage to the next. All other probabilities in Table 3b refer to a change of classification when making the transition.

For both application groups (doctoral and post-doctoral) it is striking in Table 3b that, both in going from the first (external reviewer) to the second (staff members) evaluation stage, and from the second to the third (Board of Trustees), the probability of changing from the second or third latent class to the first latent class ("award") is 0 ($\tau = 0.00$). This means an applicant only has a chance of receiving a grant from the B.I.F. if he was recommended for support at the first evaluation stage (external reviewer). Improvement to the first latent class when proceeding to the second or third evaluation stage is wholly unlikely.

The probabilities in Table 3b also indicate, however, that even applications in the first latent class in the first evaluation stage have a higher risk thereafter of not being supported by the B.I.F. in the end. True, the probability, when going on to the second evaluation stage, of remaining in the first latent class is still 52% ($\tau_{11} = 0.52$; for post-doctoral applicants) and 64% ($\tau_{11} = 0.64$; for doctoral applicants); however, when making the transition to the third stage this probability is reduced to only 18% ($\tau_{11} = 0.18$; for doctoral applicants) and 13% ($\tau_{11} = 0.13$; for post-doctoral applicants). At the same time, with every transition to the next stage, the probability of changing to the third latent class is greater, regardless of which latent class the application was in before. Even when moving to the third evaluation stage, both for doctoral ($\tau_{13} = 0.20$) and for post-doctoral applications ($\tau_{13} = 0.27$) the probability is high of changing from the first to the third latent class.

## 3.4 The influence of applicant's gender on the B.I.F. peer review process

Nonscientific statuses of the applicants (e.g., gender) are functionally irrelevant for the progress of science, and to the extent that they are used as explicit or hidden criteria in the evaluation of scientific work, the principle of universalism is being abridged (Cole, 1992, p. 162). In the last few years a series of studies have been made analyzing the influence of nonscientific statuses of applicants on judgements in the peer review (see an overview in Wood & Wessely, 2003). The status which has received the greatest attention in this research is gender. In the most influential study so far Wennerås & Wold (1997) were able to show the clear influence of the applicant's gender on the award of post-doctoral fellowships by the Swedish Medical Research Council (Stockholm).

*Comparison of two models*

To test the influence of the applicant's gender on the ratings in the B.I.F. peer review process, we calculated two latent Markov models separately for doctoral and post-doctoral fellowship applications: while the first model ($M_1$) does not differentiate between male and female applicants, the second model ($M_2$) allows for this difference. If the LR differences between $M_1$ and $M_2$ deviate significantly from zero, then there would be evidence for an influence of gender. Table 4 (above) shows for applications for doctoral fellowships ($p_{LR}^{boot}$) that the LR-differences ($M_1 - M_2$) tested by $\chi^2$-statistics are statistically significant ($M_1 - M_2$: $\Delta LR = 76.65 - 50.96 = 25.69$ $p < 0.01$, $df = 39 - 31 = 8$). The result makes it apparent that applicant's gender influences the outcome of the B.I.F. peer review process. For applications for post-doctoral fellowships on the other hand (Table 4, below) the null hypothesis for differences between $M_1$ and $M_2$ cannot be rejected. The LR does not drop statistically significantly from 40.75 ($M_1$) to 36.38 ($M_2$) ($\Delta LR = 4.37$ $p > 0.01$, $df = 6$). Hence, the outcome of the process for this application group does not appear to be significantly influenced by gender.

*Latent transition probabilities*

As the comparative model has shown a significant influence of the gender on the B.I.F. peer review process for doctoral, but *not* for post-doctoral, applicants, Table 5 shows the latent transition probabilities ($M_2$) of ratings for male and female *doctoral* applicants only (the response probabilities are omitted, because the values are quite similar to the values in Table 3b). For the transition ($t_1 - t_2$) from the first evaluation stage (external reviewer) to the second evaluation stage (staff member) applications of males have a greater chance ($\tau_{11} = 0.67$) than females to stay in the first latent class ("award"; $\tau_{11} = 0.59$). Accordingly, applications of females that are in the first latent class at the first evaluation stage change more frequently to the second ($\tau_{12} = 0.30$) or third ($\tau_{13} = 0.11$) latent class than applications of males ($\tau_{12} = 0.23$; $\tau_{13} = 0.10$). Similar differences in transition probabilities for males and females can be seen for the transition from the second to the third evaluation stage of the B.I.F. peer review process ($t_2 - t_3$). Accordingly, males have a greater chance of remaining in the first latent class ("award") - and so of being supported in the end - when progressing to the second and third evaluation stages than females. As with females, however, males have little chance of support if they have not already gained the "award" rating at the first evaluation stage.

However, the transition probabilities in Table 5 for the transition to the second evaluation stage also indicate a countervailing gender influence in the assessment of B.I.F. applications. For applications of females ($\tau_{32} = 0.36$) there is a greater probability in this transition of changing from the third to the second latent class than for applications of males ($\tau_{32} = 0.27$). Correspondingly, applications of males ($\tau_{33} = 0.73$) that at the first evaluation stage are in the third latent class are more frequently in the third latent class at the second evaluation stage than applications of females ($\tau_{33} = 0.63$). An improvement in ratings from "no award" to "possible award" in the transition from the first to the second evaluation stage is more probable for female applicants than for male applicants.

## 4   Discussion

For empirical data that can be used to check the reliability and fairness of peer review processes we can normally apply a complex structure of dependencies which have not really been appropriately statistically modelled in the past. We would highlight here the work of Jayasinghe et al. (2003), which, with its multilevel, cross-classification approach to the data analysis, first took into account that in peer review processes (1) assessors are nested into applications and (2) many assessors usually evaluated more than one application for a foundation. In this study we have applied other forms of dependencies in the peer review processes, which should be taken into account when analyzing data for testing this process.

Traditionally, in the peer review system applications go through various stages of internal and external evaluations, which are undertaken in correlation with each other. According to the Office of Management and Budget (2004, p. 12), these expensive peer reviews are appropriate for today's highly complex and multidisciplinary scientific contributions, especially those that are novel or precedent-setting. This paper presents a general methodical framework for analyzing such expensive processes by using latent Markov models. Applications for doctoral and post-doctoral fellowships that were judged in the peer review process of the B.I.F. serve as data. The B.I.F. process consists of three evaluation stages: 1) evaluation by an external reviewer, 2) internal evaluation by a staff member, and 3) final decision by the B.I.F. Board of Trustees. Using latent Markov models we examined (1) a kind of test-retest reliability of the B.I.F. peer review process. In addition, we analyzed initial latent class proportions and latent transition probabilities with the aim, (2) of suggesting ways of achieving a leaner peer review process, as well as (3) determining the influence of applicant's gender on the B.I.F. peer review process.

1. The routine cases of determining the reliability of reviewer judgements in peer review involve two or more external reviewers who judge *independently* the same scientific contribution (Cicchetti, 1991; von Eye & Mun, 2005). The weighted and unweighted Cohen's kappa and the intraclass correlation are normally used to determine this reliability (see, e.g., Daniel, 1993, 2004). However, as in multi-stage peer review processes one assessor for one evaluation stage knows the assessor's rating for the preceding evaluation stage (*dependent* ratings), these statistical measures are not appropriate. Therefore, we have determined the reliability of the peer review process by the *true* proportion of those applications for which the *dependent* ratings on the same contribution do not change from the first to the second and third stage. This proportion was based on the proportions of measurement error in the change of ratings (19% for doctoral and 18% for post-doctoral fellowship applications). For both application groups we obtained a value of ~0.80 (1 - 0.20) for the reliability of the ratings. In comparison with psychological test inventories that usually have reliability coefficients between 0.80 and 0.90, the reliability in this case can be considered satisfactory. That means applications for doctoral and post-doctoral fellowships are assessed sufficiently reliably by (1) the external reviewers, (2) the staff members, and (3) the Board of Trustees.

2. The transition probabilities of the latent Markov models for the peer review process of the B.I.F. tell us about the probabilities with which applications remain in one latent class (judgement category) in the transition from one evaluation stage to another or change to a different latent class (judgement category). The results of these transition probabilities show, both for doctoral and post-doctoral applications, that only those applicants recommended for a fellowship at the first evaluation stage (external reviewer) can obtain grant aid from the B.I.F. An improvement in the rating from "possible award" or "no award" to "award" in the transition from the first to the second and from the second to the third evaluation stage is virtually impossible.

Therefore, a form of triage or pre-screening seems desirable in the B.I.F. peer review process "in which not all grants receive the full process and deliberations of the full committee, but are rejected at an earlier stage" (Wood & Wessely, 2003, p. 32). "The goal is to allow peer reviewers to spend more time on top proposals and less effort reviewing – and re-reviewing – grants that are unlikely ever to get funded and to make reviewing a more satisfying experience" (Marshall, 1994, pp. 1212-1213). According to Marshall (1994), applicants who are rejected using triage get the message "that this is not an application that can be moved into the fundable category simply by responding to a series of complaints" (p. 1213).

Triage has been used at the NIH since 1988, after a pilot study of reviewers had suggested that they are in favour of triage (Marshall, 1994). The study of Vener et al. (1993) with empirical data from the National Cancer Institute (NCI, Bethesda, MD, USA) shows that "the conservative model [of the NIH] is valid such that the likelihood of eliminating a highly competitive application from consideration for funding is remotely small. With the model, the process of triage is fair to applicants on the one hand and is also effective in reducing consultant workloads on the other" (p. 1312).

3. The findings on the effect of gender on the B.I.F. peer review process are heterogeneous. The comparison of two latent Markov models ($M_1$ does not differentiate between male and female applicants; $M_2$ differentiates between both applicant groups) shows that applicants' gender has a statistically significant influence on the outcome of the B.I.F. peer review process in applications for doctoral fellowships. No statistically significant influence could be found for post-doctoral fellowship applications. This incongruent result reflects the heterogeneous findings of other empirical studies investigating gender effects on peer review processes. For example, some studies indicate that women scientists are at a disadvantage (Brouns, 2000; Wennerås & Wold, 1997). However, a similar number of studies report only moderate effects or no gender effects (Cole, 1992; National Science Foundation, 2000; Ward & Donnelly, 1998). An experimental study by Sonnert (1995, p. 47) found that grant submissions by women biologists received even better average evaluations than men's grant submissions did (mean rating: 3.67 vs. 3.27; $p = 0.0496$).

As, in our study, the comparative model only showed a significant influence of gender in the decision-making processes for doctoral applicants, we have only analyzed the latent transition probabilities for this group. The findings show that, although both males and females have little chance of funding if they have not already obtained the "award" rating from the external reviewer at the first evaluation stage, males however have a distinctly greater chance than females of remaining in the "award" category (and being funded) in the transitions from the first (external reviewer) to the second (staff member) and from the second to the third (Board of Trustees) evaluation stage. The transition probabilities show, however, as well as this, a contrary gender effect on the assessment of B.I.F. applications: an improvement in the rating from "no award" – given by the external reviewer at the first evaluation stage – to "possible award" (staff member) is more likely for female applicants than male applicants.

Nieva & Gutek (1980) also report similar contradictory gender effects in a research summary on the evaluation of the qualifications and performance of men and women in the work environment. "Bias … appears to work in both directions. Competent males are rated more positively than equally competent females, while incompetent males are rated lower than equally incompetent females. This pattern of results can be reconciled by the notion of sex-role congruence. Because success at most demanding situations or occupations is generally expected of males and not of females, unsuccessful females are not as heavily penalized as unsuccessful males, from whom more is expected; however, females are not rewarded for success in the same way that males are" (p. 273).

As these findings for the B.I.F. peer review process show, the suggested Markov approach can be considered as a good framework for the analysis of the peer review process. The following conditions should be met, however, in estimating Markov models for a peer review system: (1) at each evaluation stage the same rating scale is used by all assessors (internal and external) or the rating categories actually applied can be transformed for the data analysis into variables with the same categories; (2) the rating of one assessor at one evaluation stage is dependent on the assessor's rating at the preceding stage (i.e. it is known to him); and (3) the (internal and external) assessors at each evaluation stage assess the applications against the same assessment criteria (e.g., the scientific quality as demonstrated by the applicant's achievements to date, the originality of the proposed research project and the scientific standing of the laboratory where the research will be conducted).

# Footnote

**Footnote 1**

We would like to thank Dr. Rolf Langeheine, retired professor of the IPN – Leibniz Institute for Science Education Kiel (Germany) for his helpful comments on estimating the latent Markov models.

Table 1.

Contingency table of the data for the B.I.F. peer review process ($n = 1954$)

| Applications for doctoral fellowships | | | | Applications for post-doctoral fellowships | | | |
|---|---|---|---|---|---|---|---|
| External reviewer | Staff member | Board of Trustees | Observed frequencies | External reviewer | Staff member | Board of Trustees | Observed frequencies |
| 1 | 1 | 1 | 143 | 1 | 1 | 1 | 31 |
| 1 | 1 | 2 | 254 | 1 | 1 | 2 | 62 |
| 1 | 1 | 3 | 142 | 1 | 1 | 3 | 46 |
| 1 | 2 | 1 | 9 | 1 | 2 | 1 | 3 |
| 1 | 2 | 2 | 74 | 1 | 2 | 2 | 23 |
| 1 | 2 | 3 | 155 | 1 | 2 | 3 | 48 |
| 1 | 3 | 1 | 1 | 1 | 3 | 1 | 1 |
| 1 | 3 | 2 | 9 | 1 | 3 | 2 | 7 |
| 1 | 3 | 3 | 112 | 1 | 3 | 3 | 57 |
| 2 | 1 | 1 | 8 | 2 | 1 | 1 | 0 |
| 2 | 1 | 2 | 20 | 2 | 1 | 2 | 4 |
| 2 | 1 | 3 | 26 | 2 | 1 | 3 | 8 |
| 2 | 2 | 1 | 1 | 2 | 2 | 1 | 2 |
| 2 | 2 | 2 | 16 | 2 | 2 | 2 | 5 |
| 2 | 2 | 3 | 84 | 2 | 2 | 3 | 27 |
| 2 | 3 | 1 | 0 | 2 | 3 | 1 | 0 |
| 2 | 3 | 2 | 1 | 2 | 3 | 2 | 1 |
| 2 | 3 | 3 | 103 | 2 | 3 | 3 | 44 |
| 3 | 1 | 1 | 2 | 3 | 1 | 1 | 0 |
| 3 | 1 | 2 | 9 | 3 | 1 | 2 | 1 |
| 3 | 1 | 3 | 26 | 3 | 1 | 3 | 6 |
| 3 | 2 | 1 | 1 | 3 | 2 | 1 | 1 |
| 3 | 2 | 2 | 8 | 3 | 2 | 2 | 1 |
| 3 | 2 | 3 | 65 | 3 | 2 | 3 | 22 |
| 3 | 3 | 1 | 0 | 3 | 3 | 1 | 0 |
| 3 | 3 | 2 | 0 | 3 | 3 | 2 | 2 |
| 3 | 3 | 3 | 205 | 3 | 3 | 3 | 78 |

*Notes*. 1 = "award", 2 = "possible award", 3 = "no award"

Table 2.

Estimated proportions of stability and change in manifest data and latent Markov models

|  |  | Data | Markov model |
|---|---|---|---|
| Applications for doctoral fellowships | Stability | 0.24 | 0.23 |
|  | true stability |  | 0.20 |
|  | measurement error |  | 0.03 |
|  | Change | 0.76 | 0.77 |
|  | true change |  | 0.58 |
|  | measurement error |  | 0.19 |
|  | Total measurement error |  | 0.22 |
| Applications for post-doctoral fellowships | Stability | 0.22 | 0.21 |
|  | true stability |  | 0.19 |
|  | measurement error |  | 0.02 |
|  | Change | 0.78 | 0.79 |
|  | true change |  | 0.61 |
|  | measurement error |  | 0.18 |
|  | Total measurement error |  | 0.21 |

Table 3.

Estimated parameter values (and standard errors) from the latent Markov model ($n = 1954$,

75.4% applications for doctoral fellowships, 24.6% applications for post-doctoral fellowships)

a) Initial class proportion and response probabilities (row percent)

| Group | Latent class | Class proportions $\delta$'s | Response probabilities $\rho$'s | | |
|---|---|---|---|---|---|
| | | | Award | Possible award | No award |
| Applications for doctoral fellowships | 1 | 0.62 | 0.92 (0.02) | 0.06 (0.01) | 0.02 (0.01) |
| | 2 | 0.18 | 0.17 (0.03) | 0.81 (0.04) | 0.02 (0.04) |
| | 3 | 0.20 | 0.00 (n.e.) | 0.00 (n.e.) | 1.00 (n.e.) |
| Applications for post-doctoral fellowships | 1 | 0.62 | 0.92 (0.02) | 0.06 (0.01) | 0.02 (0.01) |
| | 2 | 0.18 | 0.17 (0.03) | 0.81 (0.04) | 0.02 (0.04) |
| | 3 | 0.20 | 0.00 (n.e.) | 0.00 (n.e.) | 1.00 (n.e.) |

b) Latent transition probabilities (row percent)

| Group | Latent class ($t_1$) | Latent transition probabilities $\tau$'s | | | Latent class ($t_2$) | | | |
|---|---|---|---|---|---|---|---|---|
| | | From $t_1$ to $t_2$ | | | | From $t_2$ to $t_3$ | | |
| | | Class 1 | Class 2 | Class 3 | | Class 1 | Class 2 | Class 3 |
| Applications for doctoral fellowships | 1 | 0.64 (0.02) | 0.26 (0.03) | 0.10 (0.02) | 1 | 0.18 (0.03) | 0.62 (0.05) | 0.20 (0.04) |
| | 2 | 0.00 (n.e.) | 0.54 (0.04) | 0.46 (0.04) | 2 | 0.00 (n.e.) | 0.21 (0.03) | 0.79 (0.03) |
| | 3 | 0.00 (n.e.) | 0.32 (0.03) | 0.68 (0.03) | 3 | 0.00 (n.e.) | 0.00 (n.e.) | 1.00 (n.e.) |
| Applications for post-doctoral fellowships | 1 | 0.52 (0.04) | 0.28 (0.04) | 0.20 (0.03) | 1 | 0.13 (0.05) | 0.60 (0.07) | 0.27 (0.05) |
| | 2 | 0.00 (n.e.) | 0.46 (0.07) | 0.54 (0.07) | 2 | 0.00 (n.e.) | 0.22 (0.04) | 0.78 (0.04) |
| | 3 | 0.00 (n.e.) | 0.26 (0.05) | 0.74 (0.05) | 3 | 0.00 (n.e.) | 0.04 (0.02) | 0.96 (0.02) |

*Note.* n.e. = standard error can not be calculated because of bounded parameters (0, 1)

Table 4.

Likelihood-ratio test statistic for models $M_1$ (gender differences are *not* considered in the model) and $M_2$ (gender differences are considered in the model) ($n = 1954$, 75.4% applications for doctoral fellowships, 24.6% applications for post-doctoral fellowships)

| Model | df | LR | $p_{LR}^{boot}$ |
|---|---|---|---|
| Applications for doctoral fellowships | | | |
| $M_1$ | 39 | 76.65 | 0.00 |
| $M_2$ | 31 | 50.96 | 0.02 |
| | | | |
| Applications for post-doctoral fellowships | | | |
| $M_1$ | 38 | 40.75 | 0.35 |
| $M_2$ | 32 | 36.38 | 0.26 |

*Notes.* df = degrees of freedom, LR = Log likelihood ratio, $p_{LR}^{boot}$ = bootstrapped calculated probability for loglikelihood ratio

Table 5.

Estimated latent transition probabilities (and standard errors) for male and female doctoral applicants ($n = 1474$, row percent)

| Gender | Latent class ($t_1$) | From $t_1$ to $t_2$ | | | Latent class ($t_2$) | From $t_2$ to $t_3$ | | |
| | | Class 1 | Class 2 | Class 3 | | Class 1 | Class 2 | Class 3 |
| --- | --- | --- | --- | --- | --- | --- | --- | --- |
| Male | 1 | 0.67 (0.03) | 0.23 (0.03) | 0.10 (0.02) | 1 | 0.23 (0.04) | 0.60 (0.06) | 0.17 (0.04) |
| | 2 | 0.00 (n.e.) | 0.61 (0.06) | 0.39 (0.06) | 2 | 0.00 (n.e.) | 0.21 (0.04) | 0.79 (0.04) |
| | 3 | 0.00 (n.e.) | 0.27 (0.04) | 0.73 (0.04) | 3 | 0.00 (n.e.) | 0.00 (n.e.) | 1.00 (n.e.) |
| Female | 1 | 0.59 (0.04) | 0.30 (0.04) | 0.11 (0.02) | 1 | 0.07 (0.04) | 0.67 (0.07) | 0.26 (0.05) |
| | 2 | 0.00 (n.e.) | 0.45 (0.06) | 0.55 (0.06) | 2 | 0.00 (n.e.) | 0.20 (0.04) | 0.80 (0.04) |
| | 3 | 0.01 (0.03) | 0.36 (0.05) | 0.63 (0.05) | 3 | 0.00 (n.e.) | 0.00 (n.e.) | 1.00 (n.e.) |

The heading "Latent transition probabilities $\tau$'s" spans the columns above.

*Notes*. The response probabilities and class proportions are not shown, because they are quite similar to the model for the whole data (see Table 3a).

n.e. = standard error can not be calculated because of bounded parameters (0, 1)